\documentclass[12pt]{article}
\usepackage{amssymb,amscd}
\usepackage{url}
\usepackage{stmaryrd}
\usepackage{graphicx}

\usepackage{centernot}

\newtheorem{conj}{Conjecture}[section]
\newtheorem{theo}[conj]{Theorem}
\newtheorem{rem}[conj]{Remark}

\newtheorem{defin}[conj]{Definition}
\newtheorem{prop}[conj]{Proposition}

\begin{document}
\date{August 14, 2016}
\title{\Large Symmetric  multiple chessboard complexes  \\
and a new theorem of Tverberg type}

\author{{Du\v sko Joji\'c}\\ {\small Faculty of Science}\\[-2mm] {\small University of Banja Luka}
\and Sini\v sa T. Vre\' cica\thanks{Supported by the Ministry for
Science and Technology of Serbia, Grant 174034.}
\\ {\small Faculty of Mathematics}\\[-2mm] {\small University of Belgrade}
\and Rade T. \v Zivaljevi\' c${}^\ast$\\ {\small Mathematical Institute}\\[-2mm] {\small SASA,
Belgrade}\\[-2mm]}

\maketitle 

\begin{abstract}\noindent
We prove a new theorem of Tverberg-van Kampen-Flores type
(Theorem~\ref{thm:glavna}) which confirms the conjecture of
Blagojevi\' c, Frick, and Ziegler about the existence of `balanced
Tverberg partitions' (Conjecture~6.6 in, Tverberg plus
constraints, {\em Bull.\ London Math.\ Soc.}\ 46 (2014) 953--967).
The conditions in Theorem~\ref{thm:glavna} are somewhat weaker
than in the original conjecture and we show that the theorem is
optimal in the sense that the new (weakened) condition is also
necessary. Among the consequences is a positive answer
(Theorem~\ref{thm:balanced}) to the `balanced case' of the problem
whether each {\em admissible} $r$-tuple is {\em Tverberg
prescribable}, \cite[Question~6.9]{bfz}.
\end{abstract}

\section{Introduction}

Multiple chessboard complexes are studied in their own right
\cite{krw, jvz} as interesting graph complexes where (in the
spirit of \cite{Jo-book}) the analysis of connectivity properties
is one of the central themes.

\medskip
The relevance of general chessboard complexes for theorems of
Tverberg type is well known \cite{bmz, M, VZ94, zv92, Ziegler,
Z04}.  Perhaps it should not come as a surprise, as anticipated
already in \cite{jvz}, that multiple chessboard complexes are not
an exception and that they should also quite naturally arise in
this context.

\medskip
In this paper we demonstrate that the {\em symmetric multiple
chessboard complexes} (symmetrized versions of multiple chessboard
complexes) are indeed natural configuration spaces for the proof
of new results of Tverberg type.

\medskip\noindent
Our central new results are Theorems~\ref{thm:glavna} and
\ref{thm:balanced}. Theorem~\ref{thm:glavna} provides an
affirmative answer to the following conjecture of Blagojevi\'{c},
Frick, and Ziegler, \cite{bfz}.

\begin{conj} {\em (\cite[Conjecture~6.6.]{bfz})}\label{con:bfz}
Let $r\geq 2$ be a prime power, $d \geq 1$, $N \geq (r - 1)(d +
2)$, and $r(k + 1) +s > N + 1$ for integers $k \geq 0$ and $0 \leq
s < r$. Then, for every continuous map $f : \Delta_N \rightarrow
\mathbb{R}^d$, there are $r$ pairwise disjoint faces
$\sigma_1,\ldots,\sigma_r$ of $\Delta_N$ such that
$f(\sigma_1)\cap \cdots\cap f(\sigma_r) \neq \emptyset$, with
$\textrm{dim }\sigma_i\leq k + 1$ for $1 \leq i \leq s$ and
$\textrm{dim }\sigma_i\leq k $ for $s < i \leq r$.
\end{conj}

Theorem~\ref{thm:glavna} establishes the conjecture in full
generality. Moreover, it improves the conjecture in the sense that
the condition $r(k + 1) +s
> N + 1$ in Conjecture~\ref{con:bfz} is replaced in
Theorem~\ref{thm:glavna} by a weaker and possibly more natural
condition $rk+s \geq (r-1)d$. This condition is indeed weaker
since, by assuming both conditions from Conjecture~\ref{con:bfz},
we have the inequalities,
\begin{equation}\label{eqn:opti-1}
rk+s > N+1-r \geq (r-1)(d+2)-(r-1) = (r-1)(d+1) > (r-1)d.
\end{equation}
Observe that the condition $rk+s \geq (r-1)d$ is also necessary.
Indeed it expresses the fact that if for a generic affine map $f$
the intersection $\cap_{i=1}^r~f(\sigma_i)$ is non-empty then,
\begin{equation}\label{eqn:opti-2}
 {\rm codim}(f(\sigma_1)) +\ldots + {\rm codim}(f(\sigma_r))
= {\rm codim}(\cap_{i=1}^r~f(\sigma_i)) \leq d.
\end{equation}

\begin{theo}\label{thm:glavna}
Let $r\geq 2$ be a prime power, $d \geq 1$, $N \geq (r - 1)(d +
2)$, and $rk+s \geq (r-1)d$ for integers $k \geq 0$ and $0 \leq s
< r$. Then for every continuous map $f : \Delta_N \rightarrow
\mathbb{R}^d$, there are $r$ pairwise disjoint faces
$\sigma_1,\ldots,\sigma_r$ of $\Delta_N$ such that
$f(\sigma_1)\cap \cdots\cap f(\sigma_r) \neq \emptyset$, with
$\textrm{dim }\sigma_i\leq k + 1$ for $1 \leq i \leq s$ and
$\textrm{dim }\sigma_i\leq k $ for $s < i \leq r$.
\end{theo}

The proof of Theorem~\ref{thm:glavna} is given in
Section~\ref{sec:glavna}. It relies on the shellability of the
associated configuration space (symmetric, multiple chessboard
complex),  which is established in Section~\ref{sec:shellability}.
As a consequence we obtain a fundamental connectivity result for
these complexes (Theorem~\ref{thm:treca} in
Section~\ref{sec:connectivity}). The introduction of `symmetric
multiple chessboard complexes' (Section~\ref{sec:symm}) is one of
key new ideas used in the proof of Theorem~\ref{thm:glavna}. For
this reason a more general construction of `symmetrized deleted
joins' is reviewed in Section~\ref{subsec:symm-del-joins}.

\subsection{Some consequences and relatives of
Theorem~\ref{thm:glavna}}

Theorem~\ref{thm:glavna} is apparently not an immediate
consequence of any other known result of Tverberg type (see
Section~\ref{subsec:techniques} for a brief discussion). This
observation opens the question of the relative strength of
Theorem~\ref{thm:glavna} and serves as a motivation for the
systematic study of the associated reduction methods based on
Sarkaria's inequality (Sections~\ref{sec:Constraint-Sarkaria} and
\ref{sec:A-P-problem}).

\bigskip
Theorem~\ref{thm:glavna} is a candidate for the `to date' most
general known result of van Kampen-Flores type. Here is a list of
some of its immediate consequences.

\begin{enumerate}

  \item[(1)]  It provides an optimal positive answer
(see Theorem~\ref{thm:balanced} in Section~\ref{sec:A-P-problem})
to the `balanced case' of the problem whether each {\em
admissible} $r$-tuple is {\em Tverberg prescribable}
\cite[Question~6.9]{bfz};
  \item[(2)]  The classical van Kampen-Flores theorem (see \cite[Theorem~5.1.1]{M} or \cite[Theorem~6.1]{bfz})
  is obtained if $d$ is even, $r=2$, $s=0$, and $k = \frac{d}{2}$;
  \item[(3)]  The sharpened van Kampen-Flores theorem \cite[Theorem~6.8]{bfz} corresponds to the case when
  $d$ is odd, $r=2$, $s=1$, and $k = \lfloor \frac{d}{2}\rfloor$;
  \item[(4)] The case $d=3$ of the `sharpened van Kampen-Flores theorem' is following \cite[Section~6]{bfz}
  equivalent to the Conway-Gordon-Sachs theorem which says that the complete
  graph $K_6$ on $6$ vertices is `intrinsically linked';
  \item[(5)]  The generalized van Kampen-Flores theorem \cite[Theorem~6.3]{bfz}, which improves upon
              the results of Sarkaria \cite[Section~1.5.]{sar} and Volovikov \cite{Vol-2},  follows for $s=0$
              and $k = \lceil \frac{r-1}{r}d \rceil$.
\end{enumerate}
We have already seen that the condition $rk+s\geq (r-1)d$ in
Theorem~\ref{thm:glavna} cannot be improved (see inequalities
(\ref{eqn:opti-1}) and (\ref{eqn:opti-2})). The other condition
$N\geq (r-1)(d+2)$ (in the case $k+1< d$) is also tight,  as shown
essentially by the example  \cite[Section~1.5.]{sar}, see also
\cite[Section~6]{bfz} for a related discussion.

\smallskip
Theorem~\ref{thm:glavna} can be extended to maps from the simplex
$\Delta_N$ to $d$-manifolds, following the scheme proposed by
Volovikov in \cite{Vol} (see also \cite[Section~6]{bfz}). This
extension of Theorem~\ref{thm:glavna} does not require new ideas
so we omit the details.

\subsection{Proof methods and techniques}\label{subsec:techniques}

There have been two main general proof techniques used in
topological (nonlinear) problems of Tverberg-van Kampen-Flores
type. These two proof schemes may be informally referred to as the
`direct' and `indirect' proof methods.

\begin{enumerate}
 \item[(1)] The `direct methods' rely on a variant of
{\em equivariant obstruction theory} and can be classified as:
  \begin{enumerate}
   \item[(1a)] the methods which use the high connectivity of
   the configuration space;
   \item[(1b)] the methods involving a direct calculation of the obstruction.
   \end{enumerate}
 \item[(2)] The `indirect methods' comprise two basic form
 of reductions:
   \begin{enumerate}
    \item[(2a)] the `constraint method' or the Gromov-Blagojevi\'{c}-Frick-Ziegler
    reduction;
    \item[(2b)] the methods based on `Sarkaria's index inequality'
    and its relatives.
   \end{enumerate}
  \end{enumerate}

Examples of the `direct approach' include the proof of the
`Topological Tverberg theorem' \cite{BSS} and the `Type B colored
Tverberg theorem' \cite{VZ94} (both results are classified as
(1a)). The results illustrating (1b) proof scheme are the
`Generalized van Kampen-Flores theorem' \cite{Vol-2}, the `Optimal
colored Tverberg theorem' \cite{bmz}, \cite{vz11}, and (again) the
`Topological Tverberg theorem' \cite{oz87}.

\medskip
The proof of Theorem~\ref{thm:glavna} is also direct (it is
classified as the (1a) type). It relies on the shellability of the
associated {\em symmetric multiple chessboard complex} and the
well known fact that pure $n$-dimensional shellable complexes are
$(n-1)$-connected.

\medskip The elegant and powerful `constraint method' \cite{bfz}, known also as the
Gromov-Blagojevi\'{c}-Frick-Ziegler reduction, was instrumental in
the recent spectacular construction of counterexamples for the
general topological Tverberg conjecture, see \cite{bbz} for the
exposition of the history of the problem. According to the diagram
of implications on page 2 in \cite{bz}, all known (topological)
results of Tverberg-van Kampen-Flores type are reducible by the
`constraint method' to either the {\em Topological Tverberg
Theorem} (TTT) or the {\em Optimal Colored Tverberg Theorem}
(OCTT).

\medskip
Theorem~\ref{thm:glavna} is apparently the first result which is
not an immediate consequence of either (TTT) or (OCTT).

\medskip\noindent For this reason it may be interesting to (1) explore
the nature of this phenomenon and (2) determine which results can
be reduced to Theorem~\ref{thm:glavna}, by either the `constraint
method' or some other reduction procedure. With these objectives
in mind we develop in Section~\ref{sec:Constraint-Sarkaria} some
new type (2b) reductions (based on Sarkaria's index inequality)
which appear to be better adapted to the use of `symmetrized
deleted joins' than the original `constraint method'.

\subsection{Admissible, prescribable and persistent partitions}

Theorem~\ref{thm:glavna} is an example of a Tverberg type result
where one can prescribe in advance (as in the classical van Kampen
- Flores theorem) an upper bound ${\rm dim}(\sigma_i)\leq d_i$ on
the dimension of simplices in a Tverberg $r$-partition for any map
$f : \Delta_N \rightarrow \mathbb{R}^d$.

\medskip
There is a general and a very interesting problem (Question~6.9.
in \cite{bfz}) asking for a characterization of $r$-tuples
$(d_i)_{i=1}^r$ which are `Tverberg prescribable' (see
\cite[Definition~6.7]{bfz} or our Definition~\ref{def:Tv-A-P}).

\begin{figure}[hbt]
\centering
\includegraphics[scale=0.55]{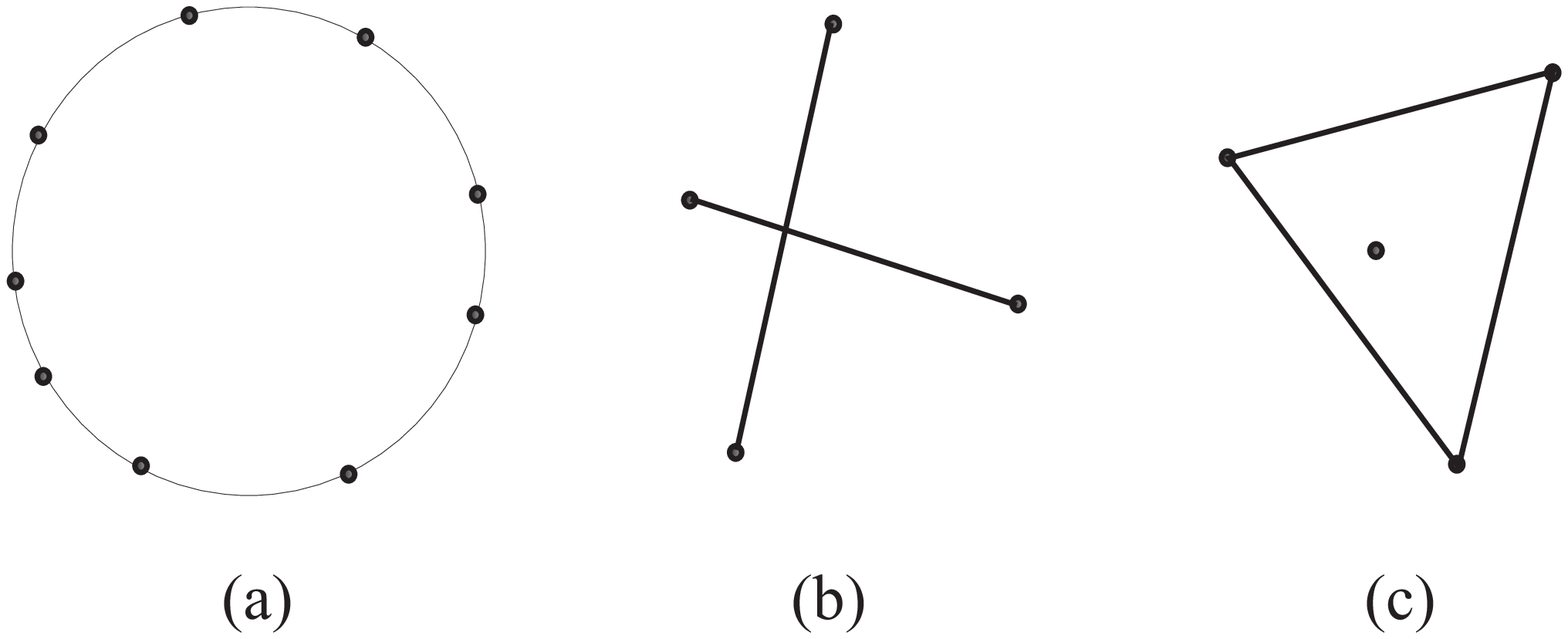}
\caption{In Radon's theorem the $(2,2)$-partitions are persistent,
while $(3,1)$ are not.} \label{fig:konveks-1}
\end{figure}

\medskip
The essence of the problem is nicely illustrated already by the
planar case of Radon's theorem, Figure~\ref{fig:konveks-1}.
Indeed, each collection $S\subset \mathbb{R}^2$ of four points in
the plane admits a partition $S = S_1\uplus S_2$ such that ${\rm
conv}(S_1)\cap {\rm conv}(S_2)\neq\emptyset$.

\smallskip
Depending on the size of the sets $S_i$ there are two types of
partitions, a $(2,2)$-partition, as in
Figure~\ref{fig:konveks-1}~(b), or a $(3,1)$-partition, as in
Figure~\ref{fig:konveks-1}~(c). Note however that these two types
of partitions are not necessarily equally represented if $S$ has
$5$ points or more. Indeed, the $(2,2)$-configurations are always
present in such a set (we may also say that $(2,2)$-partitions are
persistent) while $(3,1)$-configuration are not persistent, for
example they never appear if $S$ is in convex position,
Figure~\ref{fig:konveks-1} (a). In the language of \cite{bfz}
(Definition~6.7) the dimension pair $(d_1,d_2) = (1,1)$
(corresponding to $(2,2)$-partitions) is {\em Tverberg
prescribable} while the pair $(d_1,d_2) = (2,0)$ is not!

\bigskip
As a consequence of Theorem~\ref{thm:glavna} we establish a result
(Theorem~\ref{thm:balanced}) which addresses the case where the
prescribed bounds $d_i$ are `balanced' in the sense that $\vert
d_i - d_j \vert\leq 1$ for each $i$ and $j$.

\medskip
Theorem~\ref{thm:balanced}, as an instance of the general A-P
conjecture (Section~\ref{sec:A-P-problem}), is certainly of some
independent interest so we include in
Section~\ref{sec:A-P-problem} an alternative approach to this
result, illustrating the `extended Sarkaria's reduction'
(Section~\ref{subsec:ext-Sarkaria}).

\section{Symmetric configuration spaces}

A large symmetry group of the configuration space is one of the
key ingredients in the usual {\em configuration space/test
map}-scheme \cite{Z04} for applying topology to problems of
discrete geometry. If our first choice for the configuration space
is not sufficiently symmetric, it is quite natural to introduce a
bigger space which admits a (fixed point) free group action.

\subsection{Symmetric multiple chessboard
complexes}\label{sec:symm}

Suppose that $\mathbf{k}=(k_i)_{i=1}^n$ and
$\mathbf{p}=(p_j)_{j=1}^m$ are two sequences of non-negative
integers. The associated {\em multiple chessboard complex}
\begin{equation}\label{eqn:defin-hv}
\Delta_{m,n}^{\mathbf{k}, \mathbf{p}} =
\Delta_{m,n}^{k_1,...,k_n;p_1,...,p_m}
\end{equation}
is following \cite{krw,jvz} defined as the complex of all
rook-placements $A\subset [m]\times [n]$ such that at most $k_i$
rooks are allowed to be in the $i$-th row (for $i=1,...,n$), and
at most $p_j$ rooks are allowed in the $j$-th column (for
$j=1,...,m$). As in \cite{jvz} we pay special attention to the
complexes $\Delta_{m,n}^{\mathbf{k},\mathbf{1}} =
\Delta_{m,n}^{k_1,...,k_n; \mathbf{1}}$ where $p_1=\ldots = p_m =
1$.

\medskip
Let $G\subset S_n$ be a subgroup of the symmetric group $S_n$
acting on the set of rows of the chessboard $[m]\times [n]$. The
multiple chessboard complex $\Delta_{m,n}^{k_1,...,k_n;
\mathbf{1}}$ is rarely $G$-invariant. Since the $G$-invariance of
the configuration space is an essential feature of the usual
`configuration space/test map - scheme' \cite{Ziv-96-98, Z04}, it
is quite natural to define a symmetric version of
$\Delta_{m,n}^{k_1,...,k_n; \mathbf{1}}$.

\begin{defin}\label{def:sym-multi}
The $G$-symmetric multiple chessboard complex
\begin{equation}\label{eqn:sym-multi}
\Sigma(\Delta_{m,n}^{\mathbf{k}, \mathbf{1}}; G) =
\Sigma(\Delta_{m,n}^{k_1,...,k_n; \mathbf{1}}; G) = \bigcup_{g\in
G}~ \Delta_{m,n}^{k_{g(1)},...,k_{g(n)}; \mathbf{1}}
\end{equation}
is obtained from the multiple chessboard complex
$\Delta_{m,n}^{k_1,...,k_n; \mathbf{1}}$ by the symmetrization
with respect to $G$. In the special case $G=S_n$ we obtain the
complex,
$$\Sigma_{m,n}^{\mathbf{k}, \mathbf{1}} = \Sigma_{m,n}^{k_1,...,k_n;
\mathbf{1}} = \Sigma(\Delta_{m,n}^{k_1,...,k_n; \mathbf{1}}) :=
\Sigma(\Delta_{m,n}^{k_1,...,k_n; \mathbf{1}}; S_n)
$$
which is referred to as the {\em symmetric multiple chessboard
complex}.
\end{defin}

\begin{rem}{\rm
If one starts with more general multiple chessboard complexes
$\Delta_{m,n}^{\mathbf{k}, \mathbf{p}} =
\Delta_{m,n}^{k_1,...,k_n; p_1,..., p_m}$ it may be sometimes more
natural to symmetrize with respect to the group $G = H\times Q$
where $H$ permutes the rows and $Q$ permutes the columns of the
chessboard $[m]\times [n]$. In this paper we do not need these
more general complexes.}
\end{rem}

\subsection{Symmetrized deleted joins}\label{subsec:symm-del-joins}

The `symmetric multiple chessboard complex'
$\Sigma(\Delta_{m,n}^{k_1,...,k_n; \mathbf{1}})$, described in
Definition~\ref{def:sym-multi}, exemplifies a much more general
`symmetrized deleted join' construction.

\begin{defin}
Let $\mathcal{K} = \langle K_1,\dots, K_n\rangle$ be a collection
of not necessarily distinct simplicial complexes $K_i\subset
2^{[m]}=\Delta([m])$. Let $\mathcal{K}^{\ast n}_\Delta =
K_1\ast_\Delta\dots\ast_\Delta K_n\subset (\Delta([m])^{\ast
n}_\Delta)\cong [n]^{\ast m}$ be the associated deleted join. The
complex obtained from $\mathcal{K}^{\ast n}_\Delta$ by the
$S_n$-symmetrization,
\begin{equation}
\Sigma(\mathcal{K}^{\ast n}_\Delta):= \bigcup_{\pi\in
S_n}~K_{\pi(1)}\ast_\Delta\dots\ast_\Delta K_{\pi(n)} \subset
[n]^{\ast m}
\end{equation}
is referred to as the {\em symmetrized deleted join} of
$\mathcal{K}$.
\end{defin}

\begin{defin}
Let $\mathcal{K} = \langle K_1,\dots, K_n\rangle$ be an ordered
collection of simplicial subcomplexes of $2^{[m]} = \Delta([m])$.
Let $a = (A_1,\dots, A_n)$ be a pairwise disjoint family of
subsets of $[m]$, describing a simplex in $\Delta([m])^{\ast
n}_\Delta$. The associated $\in$-graph,
\begin{equation}
 \Gamma_{a}^\in = \Gamma_{a, \mathcal{K}}^\in =
 \{(i,j)\in [n]\times [n] \mid A_i\in K_j\}
\end{equation}
is referred to as the {\em `epsilon graph'} of $a$ and
$\mathcal{K}$ (or of $a$ alone if $\mathcal{K}$ is fixed and clear
from the context).
\end{defin}

The following proposition (relevant for
Sections~\ref{subsec:Sarkaria} and \ref{subsec:ext-Sarkaria})
provides a simple criterion for an $n$-tuple $a = (A_1,\dots,
A_n)\in [n]^{\ast m}$ to be an element of
$\Sigma(\mathcal{K}^{\ast n}_\Delta)$.

\begin{prop}\label{prop:eps-graph-crit}
A simplex $a = (A_1,\dots, A_n)\in [n]^{\ast m}\in [n]^{\ast m}$
belongs to $\Sigma(\mathcal{K}^{\ast n}_\Delta)$ if and only if
the associated $\in$-graph $\Gamma_a^\in$ admits a complete
(perfect) matching.
\end{prop}

The bipartite graph $\Gamma_a^{\notin} = \{(i,j)\in [n]\times [n]
\mid A_i\notin K_j\}$,  complementing $\Gamma_a^\in$ in the
complete bipartite graph $K_{n,n}$, is referred to as the
`non-$\in$ graph' of $a$. The following proposition illustrates
the role of non-$\in$ graphs in the definition of `collective
$n$-unavoidable complexes', introduced in \cite{jnpz}.

\begin{prop}
A collection $\mathcal{K} = \langle K_1,\dots, K_n\rangle$ of
simplicial subcomplexes of $2^{[m]} = \Delta([m])$ is {\em
`collective $n$-unavoidable'},  {\em \cite[Definition~3.1]{jnpz}},
if and only if for each simplex $a =(A_1,\dots, A_n)\in [n]^{\ast
m}$ the associated non-$\in$ graph $\Gamma_a^{\notin}$ does NOT
admit a complete matching.
\end{prop}

\medskip\noindent
{\bf Proof:} By \cite[Definition~3.1]{jnpz} a collection
$\mathcal{K} = \langle K_1,\dots, K_n\rangle$ is `collective
$n$-unavoidable' if for each (ordered) partition
$B_1\uplus\dots\uplus B_n = [m]$ of $[m]$ there exists $i\in [n]$
such that $B_i\in K_i$. This condition is easily seen to be
equivalent to the condition that the associated non-$\in$ graph
does not have a perfect matching. \hfill $\square$

\section{The connectivity of $\Sigma_{m,n}^{\mathbf{k}, \mathbf{1}}$}
\label{sec:connectivity}

One of the main results of \cite{jvz} is a lower bound for the
connectivity of the multiple chessboard complex
$\Delta_{m,n}^{\mathbf{k}, \mathbf{1}} =
\Delta_{m,n}^{k_1,...,k_n;\mathbf{1}}$.

\begin{theo}{\em (\cite[Theorem~3.1.]{jvz})}\label{thm:main}
\label{prva} The generalized chessboard complex
$\Delta_{m,n}^{k_1,...,k_n;\mathbf{1}}$ is $\mu$-connected where,
 \begin{equation}\label{eqn:glavna-ocena-jvz}
\mu = \min \{m-n-1,k_1+\cdots +k_n-2\}
\end{equation}
In particular if $m\geq k_1+\cdots +k_n+n-1$ then
$\Delta_{m,n}^{k_1,...,k_n;\mathbf{1}}$ is $(k_1+\cdots
+k_n-2)$-connected.
\end{theo}

\begin{rem}
\label{druga} {\rm As remarked already in \cite{jvz} the estimate
$\mu\geq m-n-1$ in (\ref{eqn:glavna-ocena-jvz}) can be
significantly improved for small values of $m$. However, as in
\cite{jvz}, we are here mainly interested in the values of $m$ for
which the complex $\Delta_{m,n}^{k_1,...,k_n;\mathbf{1}}$ is
$(k_1+\cdots +k_n-2)$-connected since this is precisely the result
used in applications to generalized Tverberg theorems. For this
reason our {\em working hypothesis} will be most of the time the
inequality,
\begin{equation}\label{eqn:working-hypo}
m\geq k_1+\cdots +k_n+n-1
\end{equation}
A new approach to this (central) case of Theorem~\ref{thm:main},
based on {\em discrete Morse theory}, is developed in \cite{jnpz}.
}
\end{rem}

Motivated by Conjecture~6.6. from \cite{bfz}
(Conjecture~\ref{con:bfz}), we focus our attention to the multiple
chessboard complex $K_1=\Delta_{m,n}^{\nu+1,...,\nu+1, \nu,...,
\nu;\mathbf{1}}$ of rook placements allowing at most $\nu +1$
rooks in the first $s$ rows, and at most $\nu$ rooks in the
remaining $n-s$ rows. By Theorem \ref{prva} the complex $K_1$ is
$\mu$-connected where $\mu = \min \{m-n-1,\nu n+s-2\}$.
Consequently $K_1$ is $\mu$-connected where $\mu = \nu n+s-2$,
provided the following equivalent of the inequality
(\ref{eqn:working-hypo}) is satisfied,
\begin{equation}\label{eqn:working-hypo-bis}
m\geq n(\nu+1) + s - 1
\end{equation}

\medskip
Let us suppose that $n=p^\alpha$ is a prime power. The abelian
group $G=\left( \mathbb{Z}/p\mathbb{Z}\right)^\alpha$ acts without
fixed points on the rows of the chessboard $[m]\times [n]$. By
symmetrization of $K_1$ with respect to the group $G$ we obtain
(Definition~\ref{def:sym-multi}) the symmetrized multiple
chessboard complex,
\begin{equation}
\Sigma_G = \Sigma(\Delta_{m,n}^{\nu+1,...,\nu+1, \nu,...,
\nu;\mathbf{1}}; G).
\end{equation}
More explicitly $\Sigma_G = K_1\cup \cdots \cup K_n$ is the union
of $n$ complexes where each $K_i$ is obtained from $K_1$ by the
corresponding permutation $g\in G$, in particular the simplices of
$\Sigma_G$ are rook placements on the chessboard $[m]\times [n]$
where up to $\nu+1$ rooks are allowed in some $s$ rows and at most
$\nu$ rooks in the remaining $n-s$ rows. Which combinations are
allowed is governed by the action of group $G$. Note that the
group $G=\left( \mathbb{Z}/p\mathbb{Z}\right)^\alpha$ acts without
fixed points on the complex $\Sigma_G$.

\smallskip
We find it more convenient to symmetrize with respect to the full
symmetry group so we focus our attention to the complex,
$$
\Sigma = \Sigma_{m,n}^{\nu+1,...,\nu+1, \nu,..., \nu;\mathbf{1}} =
\Sigma(\Delta_{m,n}^{\nu+1,...,\nu+1, \nu,..., \nu; \mathbf{1}};
S_n).
$$
The action of the group $G=\left(
\mathbb{Z}/p\mathbb{Z}\right)^\alpha$ on the complex $\Sigma$ is
fixed-point-free as before.

\medskip\noindent
We are interested in estimating the connectivity of the complex
$\Sigma = \Sigma_{m,n}^{\nu+1,...,\nu+1,\nu,...,\nu;\mathbf{1}}$.
Theorem~\ref{thm:treca} provides such a result and shows that
$\Sigma$ has the same connectivity lower bound as the complex $K_1
= \Delta_{m,n}^{\nu+1,...,\nu+1,\nu,...,\nu;\mathbf{1}}$.

\begin{theo}
\label{thm:treca} Suppose that,
$$
\Sigma = \Sigma_{m,n}^{k_1,...,k_s,k_{s+1},...,k_n; \mathbf{1}} =
\Sigma_{m,n}^{\nu+1,...,\nu+1,\nu,...,\nu;\mathbf{1}}
$$
is the symmetric multiple chessboard complex obtained by the
$S_n$-symmetrization of the multiple chessboard complex $K_1 =
\Delta_{m,n}^{\nu+1,...,\nu+1,\nu,...,\nu;\mathbf{1}}$ where
$k_1=\ldots = k_s = \nu+1$ and $k_{s+1}=\ldots = k_n = \nu$.
Assume that the following inequality {\em (}inequality {\em
(\ref{eqn:working-hypo-bis})}{\em )} is satisfied,
\begin{equation}
m\geq n(\nu+1) + s - 1.
\end{equation}
Then the complex $\Sigma$ is $\mu$-connected where,
\begin{equation}
\mu = k_1+\cdots +k_n-2 = \nu n+s-2
\end{equation}
\end{theo}

\medskip\noindent
{\bf Proof:} By Theorem~\ref{thm:shelling-main-new} the complex
$$\Sigma = \Sigma_{m,n}^{k_1,...,k_s,k_{s+1},...,k_n; \mathbf{1}} =
\Sigma_{m,n}^{\nu+1,...,\nu+1,\nu,...,\nu;\mathbf{1}}$$ is
shellable. Since ${\rm dim}(\Sigma) = n\nu +s -1$ the complex
$\Sigma$ is $(n\nu +s -2)$-connected and the result follows.
\hfill $\square$

\section{Shellability of $\Sigma_{m,n}^{k_1,...,k_n; \mathbf{1}}$}
\label{sec:shellability}

The following theorem was proved in \cite{jvz}.

\begin{theo}\label{thm:shelling-main-old}
For $m\geq k_1+k_2+\cdots+k_n+n-1$ the complex
$\Delta^{k_1,\ldots,k_n;\mathbf{1}}_{m,n}$ is shellable.
\end{theo}

Here we show that the shelling order for the complex
$\Delta^{k_1,\ldots,k_n;\mathbf{1}}_{m,n}$, described in
\cite[Section~4]{jvz}, can be in some cases extended to a shelling
order on the symmetrization of
$\Delta^{k_1,\ldots,k_n;\mathbf{1}}_{m,n}$ with respect to the
full symmetric group $\Sigma_n$.

\begin{theo}\label{thm:shelling-main-new}
Suppose that, $$\Sigma =
\Sigma_{m,n}^{k_1,...,k_s,k_{s+1},...,k_n; \mathbf{1}} =
\Sigma_{m,n}^{\nu+1,...,\nu+1,\nu,...,\nu;\mathbf{1}}$$ is the
symmetric multiple chessboard complex obtained by the
$S_n$-symmetrization of the multiple chessboard complex
$\Delta_{m,n}^{\nu+1,...,\nu+1,\nu,...,\nu;\mathbf{1}}$ where
$k_1=\ldots = k_s = \nu+1$ and $k_{s+1}=\ldots = k_n = \nu$.
Assume the inequality $m\geq n(\nu+1) + s - 1$ {\em ({\em the
condition}} {\em (\ref{eqn:working-hypo-bis})} in {\em {\em
Section~\ref{sec:connectivity}})}. Then the complex $\Sigma$ is
{\em shellable}.
\end{theo}

\medskip\noindent
{\bf Proof:} The dimension of the simplicial complex $\Sigma$
satisfies the inequality,
$${\rm dim}(\Sigma)+ 1 = n\nu + s \leq m - n +1.$$
By construction the complex $\Sigma$ is the union,
$$\Sigma=\bigcup\Delta^{a_1,\ldots,a_n; \mathbf{1}}_{m,n},$$ where
$(a_1,\ldots, a_n)$ is an arbitrary permutation of $(k_1,\ldots,
k_n) = (\nu+1,...,\nu+1,\nu,...,\nu)$. In other words
$a_1+a_2+\cdots+a_n = n\nu + s \leq m - n +1$, $a_i\leq \nu+1$ for
each $i$, and at most $s$ of the parameters $a_i$ can be equal
$\nu+1$. Observe that as a consequence of
Theorem~\ref{thm:shelling-main-old} all of the constituent
complexes $\Delta^{a_1,\cdots,a_r; \mathbf{1}}_{m,n}$ of $\Sigma$
are shellable.

\medskip
A facet of $\Sigma$ is naturally  encoded as an $n$-tuple
$(A_1,A_2,\ldots,A_n)$ of disjoint subsets of $[m]$, each of size
$\leq \nu+1$, such that exactly $s$ of them have $\nu+1$ elements.
We define the shelling order for the facets of $\Sigma$ by the
following construction refining the construction from \cite{jvz}.

\medskip

$F=(A_1,A_2,\ldots,A_n)$ is a predecessor of
$F'=(B_1,B_2,\ldots,B_n)$ if either,

\begin{itemize}
\item[(a)] $|A_i|=|B_i|$ for $i=1,\ldots,i_0-1$ and
$|A_{i_0}|>|B_{i_0}|$, or

\item[(b)] $|A_i|=|B_i|=a_i$ for all $i=1,2,\ldots,n$ and in the
shelling order of $\Delta^{a_1,\cdots,a_n; \mathbf{1}}_{m,n}$ the
simplex $F$ is a predecessor $F'$.

\end{itemize}
Assume that $F=(A_1,A_2,\ldots,A_n)$ precedes
$F'=(B_1,B_2,\ldots,B_n)$ in the given order. We are looking for a
facet $F''$ that is a predecessor of $F'$ such that for a vertex
$v$ of $F'$,
\begin{equation}\label{eqn:shell-1}
F\cap F'\subseteq F''\cap F'=F'\setminus\{v\}.
\end{equation}
\begin{itemize}
\item[(a)] If $|A_i|=|B_i|$ for $i=1,\ldots,i_0-1$ and
$|A_{i_0}|>|B_{i_0}|$, then there exists $j>i_0$ such that
$|B_j|>|A_j|$. Hence, there is a vertex $v$ in $B_j$ ($v$ is a
square $(x,j)$ in the $j$-th row) that does not belong to $A_j$.
Let $B'_{i_0}=B_{i_0}\cup\{x\}$ and $B'_j=B_j\setminus \{x\}$. The
facet $F''=(B_1,\ldots, B'_{i_0},\ldots,B'_j,\ldots,)$ is clearly
a predecessor of $F'$. Moreover, the facets $F,F'$ and $F''$
clearly satisfy the relation (\ref{eqn:shell-1}).

\item[(b)] If $|A_i|=|B_i|=a_i$ for all $i=1,2,\ldots,n$ then the
existence of $F''$ and $v$ follows from the shelling of the
complex $\Delta^{a_1,\ldots,a_r; \mathbf{1}}_{m,n}$.
\end{itemize}

\begin{figure}[thb]
\centering
\includegraphics[scale=0.80]{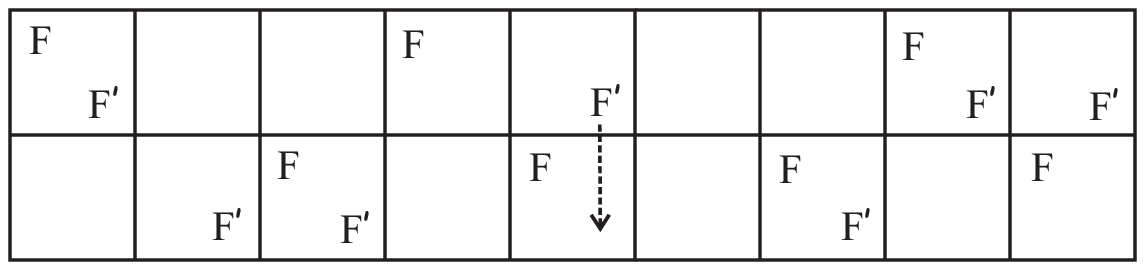}
\caption{The construction of $F''$ if $|A_i|\neq |B_i|$ for some
$i$ (the case (a)). } \label{fig:sah-sah-1}
\end{figure}

\section{Proof of Theorem~\ref{thm:glavna}}
\label{sec:glavna}

\begin{theo}\label{thm:Volovikov} {\rm (A. Volovikov \cite{Vol})}
Let $p$ be a prime number and $G = (\mathbb{Z}_p)^\alpha$ an
elementary abelian $p$-group. Suppose that $X$ and $Y$ are
fixed-point free $G$ spaces such that $\widetilde{H}^i(X,
\mathbb{Z}_p)\cong 0$ for all $i\leq n$ and $Y$ is an
$n$-dimensional cohomology sphere over $\mathbb{Z}_p$. Then there
does not exist a $G$-equivariant map $f: X\rightarrow Y$.
\end{theo}

\noindent {\bf Proof of Theorem~\ref{thm:glavna}:} We begin with
the observation that the general case of Theorem~\ref{thm:glavna},
corresponding to the pair of conditions,
\begin{equation}\label{eqn:pair-1}
N \geq (r - 1)(d + 2)  \qquad rk+s \geq (r-1)d
\end{equation}
can be reduced to the case when the second inequality is actually
an equality,
\begin{equation}\label{eqn:pair-2}
N \geq (r - 1)(d + 2)  \qquad rk+s = (r-1)d
\end{equation}
Indeed, suppose that for given $d$ and $r$ one chooses a pair
$(k', s')$ satisfying inequalities (\ref{eqn:pair-1}). Suppose
that $(k,s)$ is the unique pair satisfying the equality $rk+s =
(r-1)d$ and the condition $0\leq s < r$ ($s$ is the remainder of
$(r-1)d$ on division by $r$). Then either $k< k'$ or $k=k'$ and
$s\leq s'$. In both cases the existence of a Tverberg $r$-tuple
corresponding to the case (\ref{eqn:pair-2}) implies the existence
a Tverberg $r$-tuple corresponding to the case (\ref{eqn:pair-1})
of the theorem.

\medskip
Assuming (\ref{eqn:pair-2}) let us consider the simplicial
complex,
$$\Lambda=\Lambda(N,r;k,s)=\bigcup\Delta^{k_1,\ldots,k_r; \mathbf{1}}_{N+1,r},$$
where exactly $s$ parameters $k_i$ are equal $k+2$ and all other
parameters are $k+1$. This complex is evidently the configuration
space of all possible choices for desired Tverberg partitions with
the dimensions constrained as above. A moment's reflection reveals
that $\Lambda$ is precisely the symmetric multiple chessboard
complex $\Sigma_{m,n}^{\nu+1,...,\nu+1,\nu,...,\nu;\mathbf{1}}$
from Theorem~\ref{thm:shelling-main-new} for the following choice
of parameters,
\begin{equation}\label{eqn:parameters}
m = N+1, \quad n = r, \quad \nu = k+1.
\end{equation}

\medskip
Following the usual reduction scheme \cite[Chapter]{M} if a
desired Tverberg $r$-tuple does not exist, there exists a
$G$-equivariant map
$$F:\Lambda\rightarrow (\mathbb{R}^d)^{*r}$$
which misses the $d$-dimensional diagonal $\Delta\subset
(\mathbb{R}^d)^{*r}$, where $G = (\mathbb{Z}_p)^\alpha$ acts
without fixed points on both $\Lambda$ and
$(\mathbb{R}^d)^{*r}\setminus\Delta$. Moreover, there is a
$G$-equivariant homotopy equivalence,
$(\mathbb{R}^d)^{*r}\setminus\Delta\simeq S^D$ where $D =
(r-1)(d+1)-1$. By Theorem~\ref{thm:Volovikov} (see also
\cite[Section~6.2.6]{M} for related results of other authors),
this would lead to a contradiction if we are able to demonstrate
that the complex $\Lambda$ is $D$-connected.

\medskip
By Theorem~\ref{thm:treca} it is sufficient to guarantee the
following pair of inequalities,
\begin{equation}\label{eqn:pair-3}
m\geq n(\nu+1) + s - 1 \qquad \nu n+s-2 \geq (r-1)(d+1)-1
\end{equation}
for the choice of parameters given in (\ref{eqn:parameters}). The
second inequality in (\ref{eqn:pair-3}) is fulfilled in a tight
way as the equality $kr+s = (r-1)d$. The first inequality in
(\ref{eqn:pair-3}) is equivalent to the inequality,
\[
N+1 \geq r(k+2)+s-1 = (r-1)(d+2)+1
\]
which is precisely the first inequality in (\ref{eqn:pair-2}).
This observation completes the proof of the theorem. \hfill
$\square$

\section{Constraint method and Sarkaria's inequality}
\label{sec:Constraint-Sarkaria}

The {\em symmetric, multiple chessboard complex} played the key
role in the proof of Theorem~\ref{thm:glavna}. The apparent
incompatibility of the basic {\em `constraint method'}
(Section~\ref{subsec:Constr}) with the {\em symmetrized deleted
join construction} (Section~\ref{subsec:symm-del-joins}), may be
an initial explanation why Theorem~\ref{thm:glavna} does not
follow directly from either (TTT) or (OCTT) (see
Section~\ref{subsec:techniques}). In this section we review
alternative `indirect methods' and {\em reduction procedures},
which appear to be better adapted for the use of {\em symmetrized
deleted joins}.

\subsection{Constraint method}\label{subsec:Constr}

The proof of Theorem~\ref{thm:balanced}, based on Sarkaria's
inequality, illustrates a method which may offer an interesting
and useful alternative when the `constraint method' \cite{bfz},
(see also \cite[Section~2.9(c)]{Grom-2} and
\cite[Proposition~2.5]{L02}) is not directly applicable.

\medskip
Recall that the `constraint method' or the
Gromov-Blagojevi\'{c}-Frick-Ziegler reduction, is an elegant and
powerful method for proving results of Tverberg-van Kampen-Flores
type. It can be summarized as follows.

\begin{equation}\label{eqn:diag-constr}
\begin{CD}
K @>f>> \mathbb{R}^d\\
@VeVV @ViVV\\
L @>F>> \mathbb{R}^{d+1}
\end{CD}
\end{equation}
Suppose that $K\subset L$ are two subcomplexes of the simplex
$\Delta([n]) = 2^{[n]}$. Our objective is to prove that for each
continuous function $f : K \rightarrow \mathbb{R}^d$ there exist
disjoint faces $A_1,\dots, A_r$ in $K$ such that
$f(A_1)\cap\dots\cap f(A_r)\neq\emptyset$. To this end we assume
that $K$ is $(L,r)$-{\em unavoidable} in the sense that for each
partition $A_1\uplus\dots\uplus A_r = [n]$, where $A_j\in L$ for
each $j$, at least one of the sets $A_i$ is in $K$.

\smallskip
Assume that such a statement is known to be true for each
continuous map $F : L\rightarrow \mathbb{R}^{d+1}$. An example is
the simplex $L = \Delta^N$ (where $N = (r-1)(d+2)$ and $r = p^k$
is a prime power) by the {\em `topological Tverberg theorem'} of
B\'{a}r\'{a}ny, Shlosman, and Sz\"{ucs} \cite{BSS}, and
\"{O}zaydin \cite{oz87}  (see also \cite[Section~6.4]{M}).

\medskip
Extend the map $f : K\rightarrow \mathbb{R}^d$ to a map $\bar{f}$
($\bar{f}\circ e = f$). Suppose that $\rho : L \rightarrow
\mathbb{R}$ is the function $\rho(x) := {\rm dist}(x, K)$,
measuring the distance of the point $x\in L$ from $K$. Define $F =
(\bar{f}, \rho)  : L\rightarrow \mathbb{R}^{d+1}$ and assume that
$\Delta_1,\dots, \Delta_r$ is the associated family of vertex
disjoint faces of $\Delta^N$, such that $F(\Delta_1)\cap\ldots\cap
F(\Delta_r)\neq\emptyset$. More explicitly suppose that
$x_i\in\Delta_i$ such that $F(x_i)=F(x_j)$ for each $i,j =
1,\ldots, r$. Since $K$ is $(L,r)$-unavoidable, $\Delta_i\in K$
for some $i$. As a consequence $\rho(x_i)=0$, and in turn
$\rho(x_j)=0$ for each $j=1,\ldots, r$. If $\Delta_i'$ is the
minimal face of $\Delta^N$ containing $x_i$ then $\Delta_i'\in K$
for each $i=1,\ldots, r$ and $f(\Delta_1')\cap\ldots\cap
f(\Delta_r')\neq\emptyset$.

\medskip
For a more complete exposition and numerous examples of
applications of the `constraint method' the reader is referred to
\cite{bfz}, see also \cite[Section~2.9(c)]{Grom-2} and
\cite[Proposition~2.5]{L02}.

\subsection{The reduction based on Sarkaria's inequality}
\label{subsec:Sarkaria}

In the `direct approach' (Section~\ref{subsec:techniques}) a map
$f : K\rightarrow \mathbb{R}^d$ induces an $S_r$-equivariant map
$\hat{f} : K^{\ast r}_\Delta \rightarrow (\mathbb{R}^d)^{\ast
r}/\mathbb{R}^d$, where $\mathbb{R}^d\subset (\mathbb{R}^d)^{\ast
r}$ is the diagonal subspace. The problem is to show that this map
must have a zero, which is directly related to the problem of
estimating the equivariant index ${\rm Ind}_{G}(K^{\ast
r}_\Delta)$.\footnote{Here it is tacitly assumed, as in
\cite[Section~6]{M}, that $G=\mathbb{Z}_r$, where $r$ is a prime
number.} This transition is functorial and in particular one
obtains from (\ref{eqn:diag-constr}) the following commutative
diagram,

\begin{equation}\label{eqn:diag-del-join}
\begin{CD}
K^{\ast r}_\Delta @>\hat{f}>> (\mathbb{R}^d)^{\ast
r}/\mathbb{R}^d\\
@V\bar{e}VV @ViVV\\
L^{\ast r}_\Delta @>F>> (\mathbb{R}^D)^{\ast r}/\mathbb{R}^D
\end{CD}
\end{equation}
where $D = d+1$ (or more generally some other integer $D>d$). In
this setting the constraint method, based on the diagram
(\ref{eqn:diag-constr}), can be replaced by the index calculation
(Sarkaria's inequality) (based on the diagram
(\ref{eqn:diag-del-join})).

Recall that the Sarkaria's inequality \cite[Proposition~6.2.4]{M}
is the index inequality,
\begin{equation}\label{eqn:Sarkaria}
{\rm Ind}_G(L)\geq {\rm Ind}_G(L_0) - {\rm
Ind}_G(\Delta(L_0\setminus L)) - 1
\end{equation}
where $L_0$ is a free $G$-complex, $L\subset
 L_0$ is $G$-invariant, and $\Delta(L_0\setminus L)$ is the order
complex of the poset $(L_0\setminus L, \subseteq)$. By
specializing $L_0 = L^{\ast r}_\Delta$ and $L=K^{\ast r}_\Delta$
one obtains a tool for estimating the equivariant index of
$K^{\ast r}_\Delta$.

\smallskip
The reduction based on Sarkaria's inequality is an interesting and
useful alternative to the `constraint method', which in some cases
may facilitate calculation and lead to new observations. For
illustration here is an extension of the `van
Kampen-Flores-Gr\"{u}nbaum-Schild non-embedding theorem', proved
in \cite{jvz-2}.

\begin{theo}\label{thm:Dule-Sinisa-Rade}
Suppose that $r=p^k$ is a prime power. Let $K = K_1\ast\ldots\ast
K_s$ where each $K_i$ is an $r$-unavoidable
subcomplex\footnote{Following \cite{bfz} and \cite{jvz-2} a
complex $L\subset \Delta([k])$ is $r$-unavoidable if for each
partition $A_1\uplus\dots\uplus A_r = [k]$ of the set of vertices
there exists $i$ such that $A_i\in K$.} of the simplex
$\Delta^{m_i-1} = \Delta([m_i])$ spanned by $m_i$ vertices. Assume
that the parameters $r,d,s$ and $\{m_i\}_{i=1}^s$ satisfy the
inequality,
\begin{equation}\label{eqn:Dule-Sinisa-Rade}
(r-1)(d+s+1)+1 \leq m_1+\ldots +m_s.
\end{equation}
Then for each continuous map $f : K \rightarrow \mathbb{R}^d$
there exist $r$ disjoint faces $\Delta_1,\dots,\Delta_r$ in $K$
such that $f(\Delta_1)\cap\dots\cap f(\Delta_s)\neq\emptyset$.
\end{theo}

\subsection{Extended Sarkaria's reduction}
\label{subsec:ext-Sarkaria}

The reduction based on Sarkaria's inequality is naturally extended
to the case of symmetrized deleted joins
(Section~\ref{subsec:symm-del-joins}). This extension is
especially interesting since it is not immediately clear how to
establish a link with the `constraint method'.

\smallskip
Suppose that $\mathcal{K} = \langle K_1,\dots, K_n\rangle$ and
$\mathcal{L} = \langle L_1,\dots, L_n\rangle$ are two collections
of simplicial complexes with vertices in $[m]$ such that
$K_i\subseteq L_i$ for each $i$. In turn, there is an inclusion
map $\mathcal{K}^{\ast n}_\Delta = K_1\ast_\Delta\dots\ast_\Delta
K_n \hookrightarrow \mathcal{L}^{\ast n}_\Delta =
L_1\ast_\Delta\dots\ast_\Delta L_n$ and the inclusion $e :
\Sigma(\mathcal{K}^{\ast n}_\Delta) \hookrightarrow
\Sigma(\mathcal{L}^{\ast n}_\Delta)$ of the associated symmetrized
deleted joins (Section~\ref{subsec:symm-del-joins}).

\begin{equation}\label{eqn:symm}
\begin{CD}
\Sigma(\mathcal{K}^{\ast n}_\Delta) @>f>> (\mathbb{R}^d)^{\ast
r}/\mathbb{R}^d \\
@VeVV @ViVV\\
\Sigma(\mathcal{L}^{\ast n}_\Delta) @>F>> (\mathbb{R}^D)^{\ast
r}/\mathbb{R}^D
\end{CD}
\end{equation}

There is an associated commutative diagram (\ref{eqn:symm}) and
the equivariant index of $\Sigma(\mathcal{K}^{\ast n}_\Delta)$ can
be estimated from bellow by the inequality (\ref{eqn:Sarkaria}),
\begin{equation}\label{eqn:Sarkaria-2}
{\rm Ind}_G(\Sigma(\mathcal{K}^{\ast n}_\Delta))\geq {\rm
Ind}_G(\Sigma(\mathcal{L}^{\ast n}_\Delta)) - {\rm
Ind}_G(\Delta(Q)) - 1
\end{equation}
where $Q = \Sigma(\mathcal{L}^{\ast n}_\Delta)\setminus
\Sigma(\mathcal{K}^{\ast n}_\Delta)$.

\medskip
Propositions~\ref{prop:ind-1-nejednakost} and
\ref{prop:ind-2-nejednakost} are our first examples illustrating
Sarkaria's reduction for symmetrized deleted joins. Suppose that
$K\subset L \subset \Delta([m])=2^{[m]}$ are two simplicial
complexes on the same set of vertices. For $1\leq s\leq r$ let
$\langle K, L\rangle_s^r = \langle K,\dots, K, L,\dots, L \rangle$
be the collection $\mathcal{K}= \langle K_i\rangle_{i=1}^r$ of $r$
complexes where $K_i=L$ for $i\leq s$ and $K_j=K$ for $j>s$. Such
a collection is referred to as a $(K,L)$-collection of type
$(r,s)$. If $\mathcal{K} = \langle K_i\rangle_{i=1}^r$ is a
$(K,L)$-collection of type $(r,s)$ we denote the associated
deleted join and symmetrized deleted join by,
\[
   DJ_s^r(K,L) = \mathcal{K}^{\ast r}_\Delta \qquad
   SDJ_s^r(K,L) = \Sigma(\mathcal{K}^{\ast r}_\Delta).
\]

\begin{prop}\label{prop:ind-1-nejednakost}
Assume that $1\leq s < t <r$ where $r\geq 3$ is a prime number.
Suppose that $\langle K, L\rangle_s^r$ and $\langle K,
L\rangle_t^r$ are the associated $(K,L)$ collections of type
$(r,s)$ and $(r,t)$. Then,
\begin{equation}\label{eqn:type-1}
{\rm Ind}_G(SDJ_s^r(K,L)) \geq  {\rm Ind}_G(SDJ_t^r(K,L)) + s - t.
\end{equation}
\end{prop}

\medskip\noindent
{\bf Proof:} In light of  Sarkaria's inequality
(\ref{eqn:Sarkaria-2}) it is sufficient to show that
 \begin{equation}\label{eqn:desired}
{\rm Ind}_G(\Delta(Q))\leq t-s-1
\end{equation}
where $Q = SDJ_t^r\setminus SDJ_s^r$ and $SDJ_s^r := SDJ_s^r(K,L),
SDJ_t^r := SDJ_t^r(K,L)$. For a simplex $a = (A_1,\dots, A_r)\in
SDJ_t^r$ let $X_a:=\{i\in [r]\mid A_i\in L\setminus K\}$. Then
$a\in Q = SDJ_t^r\setminus SDJ_s^r$ if and only if $$s+1\leq \vert
X_a\vert \leq t.$$ Let $Q_s^t = \{Z\subset [r]\mid s+1\leq \vert
Z\vert\leq t\}$. Observe that ${\rm dim}(\Delta(Q_s^t))= t-s-1$
and that the action of $G = \mathbb{Z}_p$ on $Q_s^t$ is free (here
we use the fact that $t<r$).

The map $X : SDJ_t^r\setminus SDJ_s^r \rightarrow Q_s^t$ which
sends $a = (A_1,\dots, A_r)$ to $X_a\in Q_s^t$ is monotone and
$G$-equivariant. By the monotonicity and the `dimension axiom' of
the index function (see properties (i) and (v) of
\cite[Proposition~6.2.4]{M}),
\begin{equation}
 {\rm Ind}_G(\Delta(Q)) \leq {\rm Ind}_G(\Delta(Q_s^t))\leq t-s-1
\end{equation}
and the inequality (\ref{eqn:desired}) follows as an immediate
consequence. \hfill $\square$

\medskip
The condition $t<r$ limits the use of
Proposition~\ref{prop:ind-1-nejednakost} in some situations. The
following simple consequence of \cite[Proposition~6.2.4]{M}) may
be useful in such cases.

\begin{prop}\label{prop:ind-2-nejednakost}
Suppose that $\mathcal{K} = \langle K_1,\dots, K_n\rangle$ and
$\mathcal{L} = \langle L_1,\dots, L_n\rangle$ are two collections
of simplicial complexes with vertices in $[m]$ such that
$K_i\subseteq L_i$ for each $i$. If $Q = \Sigma(\mathcal{L}^{\ast
n}_\Delta)\setminus \Sigma(\mathcal{K}^{\ast n}_\Delta)$ then,
\begin{equation}\label{eqn:Sarkaria-3}
{\rm Ind}_G(\Sigma(\mathcal{K}^{\ast n}_\Delta))\geq {\rm
Ind}_G(\Sigma(\mathcal{L}^{\ast n}_\Delta)) - {\rm
dim}_G(\Delta(Q)) - 1.
\end{equation}
\end{prop}

\section{Tverberg   A-P problem}\label{sec:A-P-problem}

In this section we discuss the relation of
Theorem~\ref{thm:glavna} to the problem whether each {\em
admissible} $r$-tuple is {\em Tverberg prescribable}. This
problem, as formulated in \cite{bfz} and \cite{Bach}, will be
referred to as the Tverberg A-P problem or the Tverberg A-P
conjecture.

\begin{defin}{\em (\cite[Definition~6.7]{bfz})}\label{def:Tv-A-P}
For $d\geq 1$ and $r\geq 2$, an $r$-tuple $d = (d_1,\ldots, d_r)$
of integers is {\em admissible} if,
\begin{equation}\label{eqn:def-admissible}
\lfloor  {d}/{2}\rfloor \leq d_i\leq d   \qquad  \mbox{ {\rm and}
} \qquad  \sum_{i=1}^r (d-d_i) \leq d.
\end{equation}
An admissible $r$-tuple is {\em Tverberg prescribable} if there is
an $N$ such that for every continuous map $f : \Delta_N\rightarrow
\mathbb{R}^d$ there is a Tverberg partition $\{\sigma_1,\ldots,
\sigma_r\}$ for $f$ with $\mbox{\rm dim}(\sigma_i)=d_i$.
\end{defin}

\medskip\noindent
{\bf Question:} (Tverberg A-P problem;  \cite[Question~6.9.]{bfz})
Is every admissible $r$-tuple Tverberg prescribable?

\medskip\noindent
Note that in the Tverberg A-P problem the emphasis is on the
existence (rather then the size) of $N$ which makes a given
$r$-tuple admissible (Definition~\ref{def:Tv-A-P}).

\bigskip

The balanced case of the Tverberg A-P conjecture is the case when
the dimensions $d_1,\ldots, d_r$ satisfy the condition $\vert d_i
- d_j \vert\leq 1$ for each $i$ and $j$. In other words there
exist $0\leq s < r$ and $k$ such that $d_1=\ldots = d_s = k+1$ and
$d_{s+1}=\ldots = d_r = k$. In this case the second admissibility
condition in (\ref{eqn:def-admissible}) reduces to the inequality,
\begin{equation}\label{eqn:second-cond}
rk + s \geq (r-1)d
\end{equation}
while the first condition is redundant. The case when all
dimensions are equal $d_1=\ldots = d_r$ is answered by the
`Generalized van Kampen-Flores theorem' \cite{sar, Vol-2, bfz}.
The following theorem extends this result to the case of all
balanced $r$-tuples $d_1,\ldots, d_r$.

\begin{theo}\label{thm:balanced}
Suppose that $r=p^\kappa$ is a prime power and let $\mathbf{d} =
(d_1,\ldots, d_r)$ be a sequence of integers satisfying the
condition $\vert d_i - d_j \vert\leq 1$ for each $i$ and $j$. Then
if the sequence $\mathbf{d}$ is admissible
(Definition~\ref{def:Tv-A-P}) it is Tverberg prescribable.
\end{theo}

\medskip\noindent Theorem~\ref{thm:balanced} is an immediate consequence of
Theorem~\ref{thm:glavna}. Here we offer another proof (in the case
when $r$ is a prime) which illustrates the methods developed in
Section~\ref{subsec:ext-Sarkaria}. The second proof of the special
case of Theorem~\ref{thm:balanced} may be of some independent
interest for at least two more reasons. It shows how one can avoid
the use of Theorem~\ref{thm:glavna} by reducing it to a the
`Generalized van Kampen-Flores theorem'. Moreover, it may open a
possibility for testing some other instances of the general
Tverberg A-P conjecture.

\medskip\noindent
{\bf Proof of Theorem~\ref{thm:main}:} ($r$ is a prime number)

\noindent For $\nu = k+1$ and $1\leq s \leq r$ let $\alpha_s$ be
the $r$-sequence defined as follows,
\begin{equation}
\alpha_s(i) = \left\{\begin{array}{ll} \nu+1, & i=1,\ldots, s \\
\nu, & i=s+1,\ldots, r.
\end{array} \right.
\end{equation}
Let $\Sigma_s := \Sigma_{m,r}^{\alpha_s ; \mathbf{1}}$ be the
associated symmetric multiple chessboard complexes
(Section~\ref{sec:symm}). Observe that,
\begin{equation}\label{eqn:dim}
{\rm dim}(\Sigma_{m,r}^{\alpha_s ; \mathbf{1}}) = r\nu + s - 1 =
r(k+1)+s-1 \geq (r-1)(d+1) = {\rm dim}((\mathbb{R}^d)^{\ast
r}/\mathbb{R}^d)
\end{equation}
so in light of (\ref{eqn:symm}) and the condition
(\ref{eqn:second-cond}), it is sufficient to establish the
inequality,
 \begin{equation}\label{eqn:ind-bitno}
{\rm Ind}_G(\Sigma_{m,r}^{\alpha_s ; \mathbf{1}}) \geq r\nu + s -
1.
 \end{equation}
This is certainly true if $s = r$. Indeed, this is the case when
$\alpha_1=\dots=\alpha_r = \nu+1$ and (\ref{eqn:ind-bitno}) is
essentially the content of the `Generalized van Kampen-Flores
theorem'. Alternatively this inequality can be deduced as a direct
consequence of Theorems~\ref{thm:main} and \ref{thm:Volovikov}.

\smallskip
For $s<r$ the inequality (\ref{eqn:ind-bitno}) is established by
induction. Indeed, for $1<s\leq r$ the poset
$Q_s:=\Sigma_{m,r}^{\alpha_s ; \mathbf{1}}\setminus
\Sigma_{m,r}^{\alpha_{s-1} ; \mathbf{1}}$ is an antichain. It
follows that $\Delta(Q_s)$ is $0$-dimensional and from the
Proposition~\ref{prop:ind-2-nejednakost} we deduce that,
\begin{equation}\label{eqn:Sarkaria-4}
{\rm Ind}_G(\Sigma_{m,r}^{\alpha_{s-1} ; \mathbf{1}})\geq {\rm
Ind}_G(\Sigma_{m,r}^{\alpha_{s} ; \mathbf{1}}) - 1
\end{equation}
so if the inequality (\ref{eqn:ind-bitno}) holds for a selected
$s\leq r$ it is also true for $s-1$. \hfill$\square$

\begin{rem} {\em One can also use
Proposition~\ref{prop:ind-1-nejednakost} to deduce the case
$s=t<r-1$ of the inequality (\ref{eqn:ind-bitno}) directly  from
the case $s=r-1$. Indeed, this follows from the observation that
$\Sigma_{m,r}^{\alpha_{s} ; \mathbf{1}}$ is a $(K,L)$-collection
of the type $(r,s)$ where $K=\Delta([m])^{(k)}$ and $L =
\Delta([m])^{(k+1)}$ are the respective skeletons of the simplex
$\Delta([m]) = 2^{[m]}$. }
\end{rem}


\begin{thebibliography}{10000} {\small

\bibitem[B]{Bach} R.~Bacher et al. `An Erd\"{o}s-Szekeres-type
question', June 2011,\newline
\url{http://mathoverflow.net/questions/67762/}.

\bibitem[BBZ]{bbz} I.~B\'{a}r\'{a}ny, P.V.M.~Blagojevi{\'c}, G.M.~Ziegler.
Tverberg's theorem at 50: Extensions and counterexamples,
\textit{Notices Amer. Math. Soc.}, August 2016, 732--739.

\bibitem[BSS]{BSS}
I. B\' ar\' any, S.B. Shlosman\glossary{Shlosman, S.B.}, and
A.~Sz\H ucs.\glossary{Sz\H ucs, A.} On a topological
generalization of a theorem of Tverberg. \textit{J. London Math.
Soc.}, 23:158--164, 1981.


\bibitem[BMZ]{bmz}
P.V.M.~Blagojevi{\'c}, B.~Matschke, G.M.~Ziegler.
\newblock Optimal bounds for the colored Tverberg problem.
\newblock {\em J. European Math. Soc.}, Vol.~17, Issue 4, 2015, pp. 739–-754.

\bibitem[BFZ]{bfz}
P.V.M.~Blagojevi{\'c}, F.~Frick, G.M.~Ziegler.
\newblock Tverberg plus constraints. \newline
\newblock {\em B. London Math. Soc.}, 46:953--967, 2014.

\bibitem[BZ]{bz}
P.V.M.~Blagojevi{\'c}, G.M.~Ziegler. Beyond the Borsuk-Ulam
theorem: The topological Tverberg story, arXiv:1605.07321
[math.CO].

\bibitem[Gr10]{Grom-2} M.~Gromov. Singularities, expanders and topology of
maps. Part 2: From combinatorics to topology via algebraic
isoperimetry  {\em Geom. Funct. Anal.}  20 (2010), 416–-526.


\bibitem[JVZ-1]{jvz}
D.~Joji{\'c}, S.T.~Vre{\'c}ica, R.T.~{\v Z}ivaljevi{\'c}.
\newblock Multiple chessboard complexes and the colored Tverberg
problem,
\newblock arXiv:1412.0386 [math.CO]. \newline (accepted for publication in J. Combin. Theory Ser. A.)

\bibitem[JVZ-2]{jvz-2}
D.~Joji{\'c}, S.T.~Vre{\'c}ica, R.T.~{\v Z}ivaljevi{\'c}.
\newblock Topology and combinatorics of 'unavoidable complexes',
\newblock  arXiv:1603.08472 [math.AT].


\bibitem[JNPZ]{jnpz} D.~Joji{\'c}, I.~Nekrasov, G.~Panina, R.T.~{\v Z}ivaljevi{\'c}.
Alexander $r$-tuples  and  Bier complexes, arXiv:1607.07157
[math.CO].

\bibitem[J08]{Jo-book} J.~Jonsson. \textit{Simplicial Complexes of
Graphs}. Lecture Notes in Mathematics, Vol.\ 1928, Springer 2008.

\bibitem[KRW]{krw}
D.B. Karaguezian, V. Reiner, M.L. Wachs.
\newblock Matching Complexes, Bounded Degree Graph Complexes,
and Weight Spaces of $GL$-Complexes.
\newblock {\em Journal of Algebra} 239:77--92 2001.

\bibitem[L02]{L02} M. de Longueville. Notes on the topological Tverberg theorem,
\textit{Discrete Math.} 241 (2001) 207–-233; and \textit{Discrete
Math.} 247 (2002) 271–-297 (erratum).

\bibitem[M03]{M} J.~Matou\v sek. \textit{Using the Borsuk-Ulam Theorem.
Lectures on Topological Methods in Combinatorics and Geometry}.
Universitext, Springer-Verlag, Heidelberg, 2003. Corrected and
updated second edition,  Springer 2008.

\bibitem[\"{Oz87}]{oz87} M. \"{O}zaydin. Equivariant maps for the symmetric group.
Unpublished manuscript, available online at
\url{http://minds.wisconsin.edu/handle/1793/63829}, 1987.


\bibitem[Sar]{sar} K.S.~Sarkaria, A generalized van Kampen-Flores theorem,
\textit{Proc. Amer. Math. Soc.} 11 (1991), 559–-565.


\bibitem[Vol-1]{Vol} A. Yu. Volovikov. On a topological generalization of the
Tverberg theorem. {\em Math. Notes}, 59(3):324–-326, 1996.
Translation from {\em Mat. Zametki} 59, No.3, 454--456 (1996).

\bibitem[Vol-2]{Vol-2} A. Yu. Volovikov. On the van Kampen-Flores theorem, \textit{Math. Notes}
(5) 59 (1996), 477–-481.

\bibitem[V\v Z94]{VZ94} S.~Vre\' cica and R.~\v Zivaljevi\' c. New cases of
the colored Tverberg theorem. In H.~Barcelo and G.~Kalai, editors,
\textit{Jerusalem Combinatorics '93}, Contemporary Mathematics
Vol.\ 178, pp. 325--334, A.M.S.\ 1994.

\bibitem[V{\v{Z}}11]{vz11}
S.~Vre{\'c}ica, R.~{\v{Z}}ivaljevi{\'c}.
\newblock Chessboard complexes indomitable.
\newblock {\em J. Combin. Theory Ser. A}, 118(7):2157--2166, 2011.

\bibitem[Zi94]{zie} G.M.~Ziegler. \newblock Shellability of chessboard
complexes. {\em Israel J. Math.} 1994, Vol.~87, 97--110.

\bibitem[Zi11]{Ziegler}  G.M.~Ziegler.  \newblock $3N$ colored points in a plane.
\newblock {\em Notices of the A.M.S.} Vol.\ 58 , Number 4,
550--557, 2011.

\bibitem[{\v{Z}}V92]{zv92}
R.T. {\v{Z}}ivaljevi{\'c} and S.T. Vre{\'c}ica.
\newblock The colored {T}verberg's problem and complexes of injective
  functions.
\newblock {\em J. Combin. Theory Ser. A}, 61(2):309--318, 1992.

\bibitem[\v Ziv98]{Ziv-96-98}
R. \v Zivaljevi\' c.
\newblock User's guide to equivariant
methods in combinatorics, I and II.
\newblock {\it Publ. Inst. Math. (Beograd) (N.S.)}, (I) 59(73):114--130, 1996 and (II)
64(78):107--132, 1998.

\bibitem[\v Z04]{Z04}
R.T. \v Zivaljevi\'{c}. Topological methods. Chapter 14 in
\textit{Handbook of Discrete and Computational Geometry}, J.E.\
Goodman, J.\ O'Rourke, eds, Chapman \& Hall/CRC 2004, 305--330.

}
\end{thebibliography}
\end{document}